\documentclass[12pt]{article}
\usepackage{amsmath}\usepackage{amssymb}\usepackage{amscd}
\def\qed{\rule{6pt}{6pt}}
\def\lb{\label}
\def\be{\begin{equation}}\def\ee{\end{equation}}

\begin{document}

\begin{center}

{\LARGE{\bf Jucys-Murphy elements for Birman-Murakami-Wenzl algebras}}

\vskip .6cm

{\large {\bf A. P. Isaev}}

\vskip 0.2 cm

Bogoliubov Laboratory of Theoretical Physics, JINR\\
141980 Dubna, Moscow Region, Russia

\vskip .4cm

{\large {\bf O. V. Ogievetsky\footnote{On leave of absence from P. N. Lebedev Physical Institute,
Leninsky Pr. 53, 117924 Moscow, Russia}}}

\vskip 0.2 cm

Center of Theoretical Physics, Luminy \\
13288 Marseille, France

\vskip .4cm

\end{center}

\vspace{.5 cm}
\centerline{\bf\footnotesize Abstract}

\vspace{.6 cm}

\begin{minipage}%[10cm]
{15.4cm}

{\footnotesize The Birman-Wenzl-Murakami algebra, considered as the quotient of the braid group
algebra, possesses the commutative set of Jucys--Murphy elements. We show that the set of
Jucys--Murphy elements is maximal commutative for the generic Birman-Wenzl-Murakami algebra and
reconstruct the representation theory of the tower of Birman-Wenzl-Murakami algebras.}

\end{minipage}

\vspace{1cm}

\section{\bf Introduction}

Let $G$ be either the orthogonal group $SO(N)$ or, for $N$ even, the symplectic group $Sp(N)$. Let
${\bf g}$ be its Lie algebra and
${U}_q({\bf g})$ the corresponding quantum universal enveloping algebra. We denote the space of the defining
irreducible representation (irrep) of $G$ or ${U}_q({\bf g})$ by $V$. Let $v_i\in V$
$(i=1,\dots,N)$ be a basis of $V$. Denote by $K$ the $G$-invariant paring in $V$,
$K(v_i \otimes v_j) = K_{ij} \in {\mathbb C}$.

In 1937 R. Brauer \cite{Brau} introduced a 1-parametric family of algebras
${Br}_n(x)$ to describe the centralizer of the action of $G$ on the tensor powers $V^{\otimes n}$.
More precisely, fix the value of the parameter $x$, $x=N$. The algebra ${Br}_n(N)$ has the
representation $\tau: \, {Br}_n(N)
\rightarrow \mbox{End}(V^{\otimes n})$; the image of ${Br}_n(N)$ in this representation
coincides with the commutant of the action of $G$ on $V^{\otimes n}$. The generators of the algebra
${Br}_n(N)$ are expressed in terms
of the permutation $P$ and the operator ${\mathbf K}$ related to the $G$-invariant pairing $K$:
 $$
P(v_i \otimes v_j) = v_j \otimes v_i \; , \;\;\; {\mathbf K}(v_i \otimes v_j) = K_{ij}
\; K^{kl} \, v_k \otimes v_l \; .
 $$
 Here $K^{kl}$ is inverse to $K_{ij}$, $K^{kl}K_{lj}=\delta^k_j$.
 The Brauer algebras play the same role in the representation theory of
$SO(N)$ and $Sp(N)$  groups  as the symmetric groups in the theory of
representations of linear groups. The Brauer - Schur - Weyl duality establishes the correspondence
between the finite dimensional irreps of $SO(N)$, $Sp(N)$ and the irreps of
${Br}_n(N)$.

{}For quantum deformations ${U}_q({\bf g})$, the Brauer algebras ${Br}_n(N)$ get
$q$-deformed as well; instead of ${Br}_n(x)$ one now has a
2-parametric family of algebras $BMW_n(q,\nu)$. These algebras were introduced independently by J.
Murakami
\cite{M1} and by J. Birman and H. Wenzl \cite{BW}. The centralizers $\mbox{End}_{{U}_q({\bf
g})}(V^{\otimes n})$ are realized by specific representations $\tau$ of the Birman-Mirakami-Wenzl
algebras. The value $\nu_{_{U_q(g)}}$ of the parameter $\nu$ depends on $q$ and ${\bf g}$. In the
representation $\tau$:
$BMW_n(q,\nu_{_{U_q(g)}}) \rightarrow \mbox{End}(V^{\otimes n})$ the generators of BMW algebras
are built with the help of the Yang-Baxter $R$-operator, the $q$-analogue of the permutation $P$ 
(see \cite{Ogi} for the functorial construction of $R$-matrices of BMW-type from the 
$R$-matrices of 
$GL$-type). In contrast to the classical case, the $q$-analogue of ${\mathbf K}$ is a certain 
combination of Yang-Baxter $R$-operators, see, e.g., \cite{FRT}, \cite{Isa}.

{}For generic values of the parameters $q$ and $\nu_{_{U_q(g)}}$ the BMW algebra has the same
representations as the Brauer algebra. Different aspects of the representation theory were
extensively studied in the literature (see, e.g., \cite{IOP,W2,BB,WT,OWen}, and references
therein).

Here we generalize to the BMW algebras the approach of Vershik and Okounkov \cite{OV} developed for
the representation theory of symmetric groups and adopted to the Hecke algebra case in \cite{IsOg}.
The details of the proofs will be published in our forthcoming publication \cite{IsOg5}.

\section{\bf Birman-Murakami-Wenzl (BMW) algebras}

\subsection{Braid group and its quotients}

The braid group ${\cal B}_{M+1}$ is generated by elements $\sigma_i$, $i=1, \dots ,M$, subject to
relations:
\be
\lb{braid}
\!\!\!\!\!\!\!\!\!\!\!\!\!
{\rm Braid:} \;\;\;\; \sigma_i \, \sigma_{i+1} \, \sigma_i =
\sigma_{i+1} \, \sigma_i \,  \sigma_{i+1} \; ,
\ee
 \be
 \lb{local}
{\rm Locality:}  \;\;\; \sigma_{i} \,  \sigma_{j} = \sigma_{j} \,
\sigma_{i} \;\;\; {\mathrm{if}} \;\;  |i-j| > 1  \; .
\ee
The braid group ${\cal B}_{M+1}$ is infinite. We shall discuss certain
finite-dimensional quotients of ${\mathbb C}{\cal B}_{M+1}$. \\
{\bf 1.} The \underline{Hecke algebra $H_{M+1}(q)$} is defined by relations
\be
\lb{hecke}
(\sigma_i - q)(\sigma_i + q^{-1}) = 0 \; ,
\ee
where $q$ is a parameter; dim$(H_{M}(q)) = M!$. \\
{\bf 2.} The \underline{Birman-Murakami-Wenzl algebra $BMW_{M+1}(q,\nu)$} is defined by relations
\be
\lb{bmw01}
\left\{
\begin{array}{l}
(\sigma_i - q)(\sigma_i + q^{-1})(\sigma_i - \nu) = 0  \, , \\ [0.2cm]
\kappa_i \; \sigma_{i+1}^{\pm 1} \; \kappa_i = \nu^{\mp 1} \; \kappa_i
\;\;\; ,
\end{array}
\right.
\ee
where
\be
\lb{kappa}
\kappa_i : = \frac{(q-\sigma_i)(\sigma_i + q^{-1})}{\nu(q-q^{-1})} \; , \;\; (i = 1, \dots , M)
\; ,
\ee
$q$ and $\nu$ are parameters; dim$(BMW_{M}(q,\nu))=(2M-1)!!$.

There is a beautiful graphical presentation of the braid group and its finite dimensional
quotients. The generators $\sigma_i \in {\cal B}_{M+1}$ are depicted by

\vspace{0.2cm}

\unitlength=4.5mm
\begin{picture}(17,4)

\put(4.5,1.9){$\sigma_i \;\;  =$}

\put(7,3){$\bullet$}
\put(9,3){$\bullet$}

\put(9.7,3){$\dots$}
\put(11,3){$\bullet$}
\put(13,3){$\bullet$}
\put(16,3.1){$\dots$}
\put(15,3){$\bullet$}
\put(18,3){$\bullet$}
\put(20,3){$\bullet$}

\put(7,3.5){\tiny $1$}
\put(9,3.5){\tiny $2$}
\put(11,3.5){\tiny $i$}
\put(12.7,3.5){\tiny $i+1$}
\put(14.7,3.5){\tiny $i+2$}
\put(17.9,3.5){\tiny $M$}
\put(19.6,3.5){\tiny $M+1$}

\put(11.15,3.2){\line(1,-1){1}}
\put(12.0,2){\vector(-1,-1){0.8}}
\put(12.15,2.2){\vector(1,-1){1}}
\put(12.4,2.4){\line(1,1){0.8}}
\put(7.15,3.2){\vector(0,-1){2}}
\put(9.15,3.2){\vector(0,-1){2}}
\put(15.15,3.2){\vector(0,-1){2}}
\put(18.15,3.2){\vector(0,-1){2}}
\put(20.15,3.2){\vector(0,-1){2}}

\put(7,1){$\bullet$}
\put(9,1){$\bullet$}

\put(9.7,1){$\dots$}
\put(11,1){$\bullet$}
\put(13,1){$\bullet$}
\put(16,1.1){$\dots$}
\put(15,1){$\bullet$}
\put(18,1){$\bullet$}
\put(20,1){$\bullet$}

\put(22,1.9){.}

\end{picture}

\noindent
For the locality relation (\ref{local}) we have $(i+1 <j<M)$

\unitlength=4.5mm
\begin{picture}(17,4)

\put(1.5,1){$\sigma_i \sigma_j \;\;  =$}

\put(5.7,3){$\dots$}
\put(7,3){$\bullet$}
\put(9,3){$\bullet$}
\put(10,3){$\dots$}
\put(14.5,3){$\dots$}
\put(11.5,3){$\bullet$}
\put(13.5,3){$\bullet$}

\put(7,3.5){\tiny $i$}
\put(8.7,3.5){\tiny $i+1$}
\put(11.5,3.5){\tiny $j$}
\put(13.2,3.5){\tiny $j+1$}

\put(7.15,3.2){\line(1,-1){1}}
\put(8.0,2){\vector(-1,-1){0.8}}
\put(8.15,2.2){\vector(1,-1){1}}
\put(8.4,2.4){\line(1,1){0.8}}

\put(11.65,3.2){\vector(0,-1){2}}
\put(13.65,3.2){\vector(0,-1){2}}

\put(5.7,1){$\dots$}
\put(7,1){$\bullet$}
\put(9,1){$\bullet$}
\put(10,1.1){$\dots$}
\put(14.5,1.1){$\dots$}
\put(11.5,1){$\bullet$}
\put(13.5,1){$\bullet$}

\put(11.65,1.2){\line(1,-1){1}}
\put(12.5,0){\vector(-1,-1){0.8}}
\put(12.65,0.2){\vector(1,-1){1}}
\put(12.9,0.4){\line(1,1){0.8}}

\put(7.15,1.2){\vector(0,-1){2}}
\put(9.15,1.2){\vector(0,-1){2}}

\put(5.7,-1){$\dots$}
\put(7,-1){$\bullet$}
\put(9,-1){$\bullet$}
\put(10,-0.9){$\dots$}
\put(14.5,-0.9){$\dots$}
\put(11.5,-1){$\bullet$}
\put(13.5,-1){$\bullet$}

%%%%%%%%%%%%%%%%%%%%%%%%%%%%%%%%

\put(17,1){$  =$}

\put(20.7,3){$\dots$}
\put(22,3){$\bullet$}
\put(24,3){$\bullet$}
\put(25,3){$\dots$}
\put(29.5,3){$\dots$}
\put(26.5,3){$\bullet$}
\put(28.5,3){$\bullet$}

\put(22,3.5){\tiny $i$}
\put(23.7,3.5){\tiny $i+1$}
\put(26.5,3.5){\tiny $j$}
\put(28.2,3.5){\tiny $j+1$}

\put(26.6,3.2){\line(1,-1){1}}
\put(27.45,2){\vector(-1,-1){0.8}}
\put(27.6,2.2){\vector(1,-1){1}}
\put(27.85,2.4){\line(1,1){0.8}}

\put(22.3,3.2){\vector(0,-1){2}}
\put(24.3,3.2){\vector(0,-1){2}}

\put(20.7,1){$\dots$}
\put(22,1){$\bullet$}
\put(24,1){$\bullet$}
\put(25,1.1){$\dots$}
\put(29.5,1.1){$\dots$}
\put(26.5,1){$\bullet$}
\put(28.5,1){$\bullet$}

\put(22.3,1.2){\line(1,-1){1}}
\put(23.15,0){\vector(-1,-1){0.8}}
\put(23.3,0.2){\vector(1,-1){1}}
\put(23.55,0.4){\line(1,1){0.8}}

\put(26.65,1.2){\vector(0,-1){2}}
\put(28.7,1.2){\vector(0,-1){2}}

\put(20.7,-1){$\dots$}
\put(22,-1){$\bullet$}
\put(24,-1){$\bullet$}
\put(25,-0.9){$\dots$}
\put(29.5,-0.9){$\dots$}
\put(26.5,-1){$\bullet$}
\put(28.5,-1){$\bullet$}

\put(32,1){$   = \;\; \sigma_j \sigma_i \; $.}
\end{picture}

\vspace{1cm}

\noindent
The braid relation is

\vspace{1.2cm}

\unitlength=4.5mm
\begin{picture}(17,4)

\put(0.5,1.9){$\sigma_{i+1} \sigma_i \sigma_{i+1} \;\;    =$}

\put(8.15,5.2){\line(1,-1){1}}
\put(9.0,4){\vector(-1,-1){0.8}}
\put(9.15,4.2){\vector(1,-1){1}}
\put(9.4,4.4){\line(1,1){0.8}}

\put(6.15,5.2){\vector(0,-1){2}}

\put(6,5){$\bullet$}
\put(8,5){$\bullet$}
\put(10,5){$\bullet$}

\put(6,5.5){\tiny $i$}
\put(7.7,5.5){\tiny $i+1$}
\put(9.7,5.5){\tiny $i+2$}

%%%%%%%%%%%%%%%%%%%%%%%%

\put(6,3){$\bullet$}
\put(8,3){$\bullet$}
\put(10,3){$\bullet$}

\put(6.15,3.2){\line(1,-1){1}}
\put(7.0,2){\vector(-1,-1){0.8}}
\put(7.15,2.2){\vector(1,-1){1}}
\put(7.4,2.4){\line(1,1){0.8}}
\put(10.15,3.2){\vector(0,-1){2}}

\put(6,1){$\bullet$}
\put(8,1){$\bullet$}
\put(10,1){$\bullet$}

%%%%%%%%%%%%%%%%%%%%%%%%%%%%%%%%

\put(8.15,1.2){\line(1,-1){1}}
\put(9.0,0){\vector(-1,-1){0.8}}
\put(9.15,0.2){\vector(1,-1){1}}
\put(9.4,0.4){\line(1,1){0.8}}

\put(6.15,1.2){\vector(0,-1){2}}

\put(6,-1){$\bullet$}
\put(8,-1){$\bullet$}
\put(10,-1){$\bullet$}

%%%%%%%%%%%%%%%%%%%%%%%%%%%%%%%%%

\put(12.5,1.9){$=$}

\put(15,5.5){\tiny $i$}
\put(16.7,5.5){\tiny $i+1$}
\put(18.7,5.5){\tiny $i+2$}

\put(15,5){$\bullet$}
\put(17,5){$\bullet$}
\put(19,5){$\bullet$}

\put(15.15,5.2){\line(1,-1){1}}
\put(16.0,4){\vector(-1,-1){0.8}}
\put(16.15,4.2){\vector(1,-1){1}}
\put(16.4,4.4){\line(1,1){0.8}}
\put(19.15,5.2){\vector(0,-1){2}}

\put(15,3){$\bullet$}
\put(17,3){$\bullet$}
\put(19,3){$\bullet$}

%%%%%%%%%%%%%%%%%%%%%%%%%%%%%%%%

\put(17.15,3.2){\line(1,-1){1}}
\put(18.0,2){\vector(-1,-1){0.8}}
\put(18.15,2.2){\vector(1,-1){1}}
\put(18.4,2.4){\line(1,1){0.8}}

\put(15.15,3.2){\vector(0,-1){2}}

\put(15,1){$\bullet$}
\put(17,1){$\bullet$}
\put(19,1){$\bullet$}

%%%%%%%%%%%%%%%%%%%%%%%%%%%%%%%%%%%%%%%

\put(15,-1){$\bullet$}
\put(17,-1){$\bullet$}
\put(19,-1){$\bullet$}

\put(15.15,1.2){\line(1,-1){1}}
\put(16.0,0){\vector(-1,-1){0.8}}
\put(16.15,0.2){\vector(1,-1){1}}
\put(16.4,0.4){\line(1,1){0.8}}
\put(19.15,1.2){\vector(0,-1){2}}

\put(20.5,1.9){$  = \;\; \sigma_i \sigma_{i+1} \sigma_i \;\;$.}

\end{picture}

\vspace{1cm}

\noindent
It is sometimes convenient to depict the element (\ref{kappa}) by

\vspace{0.4cm}

\unitlength=4.5mm
\begin{picture}(17,4)

\put(4,1.9){$\kappa_i \;\;  =$}

\put(6.5,3){$\bullet$}
\put(9,3){$\bullet$}

\put(7.4,3.1){$\dots$}
\put(7.4,1.1){$\dots$}
\put(16,3.1){$\dots$}
\put(16,1.1){$\dots$}

\put(11,3){$\bullet$}
\put(13,3){$\bullet$}
\put(15,3){$\bullet$}
\put(18,3){$\bullet$}
\put(20,3){$\bullet$}

\put(6.5,3.5){\tiny $1$}
\put(8.5,3.5){\tiny $i-1$}
\put(11,3.5){\tiny $i$}
\put(12.7,3.5){\tiny $i+1$}
\put(14.7,3.5){\tiny $i+2$}
\put(17.9,3.5){\tiny $M$}
\put(19.6,3.5){\tiny $M+1$}

\put(6.65,3.2){\vector(0,-1){2}}
\put(9.15,3.2){\vector(0,-1){2}}
\put(15.15,3.2){\vector(0,-1){2}}
\put(18.15,3.2){\vector(0,-1){2}}
\put(20.15,3.2){\vector(0,-1){2}}

\put(6.5,1){$\bullet$}
\put(9,1){$\bullet$}

\put(11,1){$\bullet$}
\put(13,1){$\bullet$}
\put(15,1){$\bullet$}
\put(18,1){$\bullet$}
\put(20,1){$\bullet$}

\put(12.25,3.2){\oval(2.1,1.2)[b]}
\put(12.15,1.2){\oval(2.1,1.2)[t]}

\put(22,1.9){.}

\end{picture}

\noindent
Below we shall omit the reference to the parameters in the notation $H_{M}(q)$ and $BMW_{M}(q,\nu)$
and write simply $H_{M}$ and $BMW_{M}$.

\subsection{Affine BMW algebras $\alpha BMW_{M+1}$}

Affine Birman-Murakami-Wenzl algebras  $\alpha B\! M \! W_{M+1}$ are extensions of the algebras
$BMW_{M+1}$. The algebras  $\alpha B\! M \! W_{M+1}$ are
generated by the elements $\{\sigma_1,\dots,\sigma_M\}$ with relations (\ref{braid}),
(\ref{local}), (\ref{bmw01}) and the affine element $y_1$ which satisfies
\be
\lb{bmw02}
\begin{array}{c}
\sigma_1 \, y_1 \, \sigma_1 \, y_1 = y_1 \, \sigma_1 \, y_1 \, \sigma_1 \; , \;\;\;
[\sigma_k, \, y_1]=0 \;\;\; \mathrm{for} \;\; k > 1 \; ,
\\ [0.2cm]
\kappa_1 \; y_1 \, \sigma_1 \, y_1 \, \sigma_1 = c \, \kappa_1 =
\sigma_1 \, y_1 \, \sigma_1 \,  y_{1} \; \kappa_1  \; ,
\\ [0.2cm]
\kappa_1 \; y_1^n \; \kappa_1 = \hat{z}^{(n)} \kappa_1\ \ ,
\;\;\; n=1,2,3,\dots  \; .
\end{array}
\ee
where $c$, $\hat{z}^{(n)}$ are central elements. Initially, for the Brauer algebras, the affine 
version was introduced by M. Nazarov \cite{Nazar}.

Consider the set of affine elements
$$ y_{k+1} = \sigma_k \,
y_k \, \sigma_k \; , \;\;\; k=1,2,\dots ,M \; .
$$
The elements $y_k$ $(k=1,2,\dots,M+1)$ generate a commutative subalgebra $Y_{M+1}$ in
$\alpha B\! M \! W_{M+1}$.

\subsection{Central elements in $\alpha BMW$ algebra}

We need some information about the center of $\alpha B\!M\!W$.

\noindent
{\bf Proposition 1.} {\it The elements
$$
\hat{\cal Z} = y_1 \cdot y_2 \cdots y_{M} \; ,
\;\;\; \hat{\cal Z}^{(n)}_{M} =   \sum_{k=1}^{M} \left( y^n_{k} - c^n \,
y^{-n}_{k} \right) \;\; ,\;\; n \in {\mathbb N}\ \ ,
$$
are central in the $\alpha B\!M\!W_{M}$ algebra.}

\noindent
{\bf Remark.} The set of "power sums" $\hat{\cal Z}^{(n)} = \sum_k (y_{k}^n - c^n \, y_{k}^{-n})$
has the generating function
$$
{\cal Z}(t) = \sum_{n=1} \hat{\cal Z}^{(n)} t^{n-1}
  = \frac{d}{d t}  \log \left( \prod_{k=1} \frac{y_{k}- c \, t}{1-y_{k}t} \right) \ .
$$

Consider an ascending chain of subalgebras
$$
\alpha B\!M\!W_0\subset \alpha B\!M\!W_1\subset
\alpha B\!M\!W_2\subset\dots\subset
\alpha B\!M\!W_{M}\subset \alpha B\!M\!W_{M+1} \; ,
$$
where $\alpha B\!M\!W_{0},\alpha B\!M\!W_{1}$ and $\alpha B\!M\!W_{j}$ $(j>1)$ are generated by $\{ 
c, \hat{z}^{(n)}\}$, $\{ c, \hat{z}^{(n)}, y_1\}$ and $\{ c, \hat{z}^{(n)}, y_1,
 \sigma_1,\sigma_2,
 \dots , \sigma_{j-1} \}$, respectively.
For the corresponding commutative subalgebras we have $Y_1\subset Y_2\subset\dots\subset
Y_{M}\subset Y_{M+1}$.

\noindent
{\bf Proposition 2.} {\it Let $\hat{Z}_{k}^{(n)}$ be the central elements in the algebra $\alpha
B\!M\!W_{k}$, $\alpha B\!M\!W_{k} \subset \alpha B\!M\!W_{k+2}$, defined by the generating function
\be
\lb{gfun}
\begin{array}{l}
{\displaystyle  \sum_{n=0}^\infty \hat{Z}_{k}^{(n)} t^{n} = - \frac{\nu}{(q-q^{-1})} +
\frac{1}{(1 - c \, t^2)}  +
\left( \sum_{n=0}^{\infty} t^{n} \, \hat{z}^{(n)}  + \frac{\nu}{(q-q^{-1})} - \frac{1}{(1
-c \, t^2)} \right)  }
\\
\hspace{1.2cm}
{\displaystyle \cdot \prod_{r=1}^{k}
\frac{(1-y_{r} t)^2(q^2  - c \, y^{-1}_r t) (q^{-2}  -c \, y^{-1}_r t)}{
(1-c \, y^{-1}_r t)^2(q^2  - y_r t) (q^{-2}  - y_r t)} } \; .
\end{array}
\ee
The following relations hold
$$
\kappa_{k+1} y_{k+1}^{n} \kappa_{k+1} = \hat{Z}_{k}^{(n)} \, \kappa_{k+1} \in \alpha B\!M\!W_{k+2}
\;\;\;\; (\hat{Z}_{0}^{(n)} \equiv \hat{z}^{(n)}) \; .
$$}

\noindent
 {\bf Remark.}
The evaluation map $\alpha B\!M\!W_{M}\rightarrow B\!M\!W_{M}$ is defined by
 \be
 \lb{evmap}
  y_1\mapsto 1 \;\; \Rightarrow \;\; c \mapsto \nu^2 \; , \;\;\;
  \hat{z}^{(n)} \mapsto 1+\frac{\nu^{-1} -\nu}{q-q^{-1}} \; .
 \ee
Under this map the function (\ref{gfun}) transforms into the generating function presented in
\cite{BB}.

\subsection{Intertwining operators in $\alpha B\!M\!W$ algebra}

Introduce the {\it intertwining} elements $U_{k+1} \in  \alpha B\!M\!W_{M+1}$
$(k=1,\dots,M)$
\be
\lb{intw}
U_{k+1} = [\sigma_k , \, y_{k} - c \, y_{k+1}^{-1}]  \; .
\ee

\noindent
{\bf Proposition 3.} {\it The elements $U_k$ satisfy
$$
U_{k+1} y_k = y_{k+1} \, U_{k+1} \; , \;\;\; U_{k+1} y_{k+1} = y_{k} \, U_{k+1} \; ,
\;\;\;
U_{k+1} y_i = y_{i} \, U_{k+1} \;\;\; {\mathrm{for}}\ \  i \neq k,k+1 \; ,
$$
\be
\lb{bmw05}
U_{k+1} \, [\sigma_k , \, y_{k}] = (q y_k - q^{-1} y_{k+1})(q y_{k+1} - q^{-1} y_{k})
\left(1 - \frac{c}{y_k \, y_{k+1}} \right) \; ,
\ee
$$
U_{k+1} \, U_k \, U_{k+1} = U_k \, U_{k+1} \, U_k  \; ,
$$
$$
\kappa_k \, U_{k+1} = U_{k+1} \, \kappa_k = 0 \; .
$$}
The elements $U_{k}$ provide an important information about the spectrum of the affine elements
$\{ y_j \}$.

\noindent
{\bf Lemma 1.} {\it The spectrum of the elements $y_j \in  \alpha B\!M\!W_{M+1}$ satisfies
\begin{equation}
\label{spec1}
{\rm Spec} (y_j) \subset \{ q^{2 {\mathbb Z}} \cdot {\rm Spec} (y_1), \; c \, q^{2 {\mathbb Z}}
\cdot {\rm Spec} (y_1^{-1}) \} \; ,
\end{equation}
where ${\mathbb Z}$ is the set of integer numbers.}

\noindent
{\bf Proof.} Induction in $j$. Eq. (\ref{spec1}) obviously holds for $y_1$. Assume that
$$
{\rm Spec} (y_{j-1}) \subset \{ q^{2 {\mathbb Z}} \cdot {\rm Spec} (y_1), c \, q^{2 {\mathbb Z}}
\cdot {\rm Spec} (y_1^{-1}) \} \; , \;\; j >1 \; .
$$
Let $f$ be the characteristic polynomial of $y_{j-1}$, $f(y_{j-1}) = 0$. Then
$$
\begin{array}{c}
0=U_j f(y_{j-1}) [\sigma_{j-1} , y_{j-1} ] = f(y_{j}) U_j [\sigma_{j-1} , y_{j-1} ] \\[.5em]
= f(y_{j}) (q^2 y_{j-1} - y_{j})( y_{j} - q^{-2} y_{j-1})
\left( y_{j} - c \, y^{-1}_{j-1}  \right)  y_{j}^{-1} \; .
\end{array}
$$
Here we used (\ref{bmw05}). Thus
${\rm Spec} (y_j) \subset {\rm Spec}(y_{j-1}) \cup q^{\pm 2} \cdot {\rm Spec}(y_{j-1}) \cup
c \cdot {\rm Spec}(y_{j-1}^{-1})$. \hfill \qed

\vspace{0.2cm}

  We denote the image of $w \in \alpha B\!M\!W_{M}$ under the evaluation map (\ref{evmap}) by
$\tilde{w}$, e.g., $y_j \mapsto \tilde{y}_j$.
The Jucys-Murphy (JM) elements $\tilde{y}_j$ $(j=2,\dots,M)$ are the images of $y_j$:
$$
\tilde{y}_{j}=\sigma_{j-1}\dots\sigma_2\,\sigma_1^2\,\sigma_2\dots\sigma_{j-1}\;\in\; B\!M\!W_{M}\; .
$$
Lemma 1 provides the information about the spectrum of JM elements
$\tilde{y}$'s.

\noindent
{\bf Corollary.} Since $\tilde{y}_1=1$ and $\tilde{c}=\nu^2$, it follows from (\ref{spec1}) that
 \be
\lb{spec2}
{\rm Spec} (\tilde{y}_j) \subset \{ q^{2 {\mathbb Z}} , \; \nu^2 q^{2 {\mathbb Z}} \} \; .
 \ee

\section{\bf Representations of affine algebra $\alpha B\!M\!W_2$}

\subsection{$\alpha B\!M\!W_2$ algebra and its modules $V_D$}

The elements $\{y_i, y_{i+1}, \sigma_i, \kappa_i \} \in \alpha B\!M\!W_M$ (for fixed $i < M$)
satisfy
\be
\lb{bmw06}
(q-q^{-1}) \kappa_i = \sigma_i^{-1} - \sigma_i + (q-q^{-1}) \; ,
\ee
\be
\lb{bmw08}
y_{i+1} = \sigma_i y_i \sigma_i \; , \;\;\; y_i y_{i+1} = y_{i+1} y_i \; ,
\;\;\;
 \kappa_i y_i^n \kappa_i = \hat{Z}_{i-1}^{(n)} \kappa_i \; ,
\ee
\be
\lb{bmw07}
y_i y_{i+1} \kappa_i = c \, \kappa_i
= \kappa_i y_{i+1} y_i \; .
\ee
The elements $c$ and $\hat{Z}_{i-1}^{(n)}$ commute with
$\{y_i, y_{i+1}, \sigma_i, \kappa_i \}$. The elements $\{y_i, y_{i+1},
\sigma_i, \kappa_i \} \in \alpha B\!M\!W_M$ generate a subalgebra isomorphic to $\alpha B\!M\!W_2 $.

Below we investigate representations $\rho$ of $\alpha B\!M\!W_2$ for which the generators
$\rho(y_i)$ and
$\rho(y_{i+1})$ are diagonalizable and $\rho(c)=\nu^2 \cdot {\text{Id}}$. Let $\psi$
be a common eigenvector of $\rho(y_i)$ and
$\rho(y_{i+1})$ with some eigenvalues $a$ and $b$:
$$
\rho(y_i) \, \psi = a \, \psi \; , \;\;\; \rho(y_{i+1}) \, \psi = b
\, \psi \; .
$$
The element $\hat{z}=y_i y_{i+1}$ is central in $\alpha \! B\!M\!W_2$. There are two possibilities:
\be
\lb{bmw09}
\begin{array}{l}
{\mathbf 1.} \;\;\;\; \rho(\kappa_i) \neq 0 \;\;\; \stackrel{\rm eq. \;
(\ref{bmw07})}{\Longrightarrow}
\;\;\; \rho(y_i y_{i+1})  = \nu^2 \cdot {\text{Id}} \;\; \Rightarrow \;\; \underline{a \, b = \nu^2}; \\ [0.2cm]
{\mathbf 2.} \;\;\;\; \rho(\kappa_i) = 0 \;\; , \;\;\underline{{\mathrm{the\ product}}\ \ a \, b
\;\; {\rm is} \; {\rm not} \; {\rm fixed}}.
\end{array}
\ee

To save space we shall often omit the symbol $\rho$ and denote, slightly abusively, the operator
$\rho(x)$ for $x \in \alpha B\!M\!W$ by the same letter $x$; this should not lead to a confusion.

Applying the operators from $\alpha B\!M\!W_2$ to the vector $\psi$ we produce, in general infinite
dimensional, $\alpha B\!M\!W_2$-module $V_{\infty}$ spanned by
$$
\begin{array}{ll}
e_2 = \psi \; , & \\
e_1 = \kappa_i \psi \; , & \;\;\; e_3 = \sigma_i \psi \; , \\
e_4 = y_i \kappa_i \psi \; , & \;\;\; e_5 = \sigma_i y_i \kappa_i  \psi \; , \\
e_6 = y^2_i \kappa_i \psi \; , & \;\;\; e_7 = \sigma_i y^2_i \kappa_i  \psi \; ,  \\
\;\;\; \dots\dots\dots \; ,  & \;\;\;\;\;\;  \dots\dots\dots  \; , \\
e_{2k +2} =  y^k_i \kappa_i \psi \; , & \;\;\; e_{2k +3} =  \sigma_i y^k_i \kappa_i  \psi \;\;
(k \geq 1) \; , \dots\dots \; .
\end{array}
$$
Using relations (\ref{bmw06}) -- (\ref{bmw07}) for $\alpha B\!M\!W_2$ one can write down the left
action of elements
$\{y_i, y_{i+1}$, $\sigma_i, \kappa_i \}$ on $V_{\infty}$.
Our aim is to understand when the sequence $e_j$ can terminate giving therefore rise to a finite
dimensional module $V_D$ (of dimension  $D$) of $\alpha B\!M\!W_2$ and investigate the
(ir)reducibility of $V_D$.

We distinguish 3 cases for the module $V_D$:

 \begin{itemize}
  \item[\bf (i)]
$\kappa_i V_D=0$ (i.e., $\kappa_i \, e = 0 \;\;
\forall \; e \in V_D$) and in particular $\kappa_i \, \psi = 0$. Therefore, $e_j = 0$
for all $j\neq 2,3$ and  $V_{\infty}$ reduces to a 2-dim module with the basis $\{ e_2,e_3\}$. In
view of (\ref{bmw09}) the product $a \, b$ is not fixed and the irreps coincide with the irreps of
the affine Hecke algebra
$\alpha H_2$ considered in \cite{IsOg}.
 \item[\bf (ii)]
$\kappa_i V_D \neq 0$ (i.e., $\exists \; e \in V_D$:
$\kappa_i e \neq 0$). The module $V_D$ is extracted
from $V_{\infty}$ by constraints
 \be
 \lb{bmw10}
e_{2k+4} = \sum_{m=1}^{2k+3} \alpha_m \, e_{m} \;\;\; (k \geq -1) \; , \;\; ab = \nu^2 \; ,
\ee
with some parameters $\alpha_m$. The independent basis vectors are
$(e_1,e_2,\dots,e_{2k+3})$. The module $V_D$ has  \underline{odd dimension}.
 \item[\bf (iii)]
$\kappa_i V_D \neq 0$ and additional constraints are
 \be
 \lb{bmw11}
e_{2k+3} = \sum_{m=1}^{2k+2} \alpha_m \, e_{m} \;\;\; (k \geq 0) \; , \;\; ab = \nu^2 \; .
 \ee
The independent basis vectors are $(e_1,e_2,\dots,e_{2k+2})$. The module
$V_D$ has \underline{even dimension}.
 \end{itemize}
Below we consider a version $\alpha' BMW_2$ of the affine BMW algebra. The additional requirement
for this algebra concerns the spectrum of $y_i,y_{i+1} \in \alpha' BMW_2$:
$$
{\rm Spec} (y_j) \subset \{ q^{2 {\mathbb Z}} , \; \nu^2 q^{2 {\mathbb Z}} \} \; .
  $$
  The evaluation map (\ref{evmap}) descends to the algebra $\alpha' BMW$ (cf. Corollary 
  after Lemma 1).
  In particular for the cases {\bf (ii)} and {\bf (iii)} we have
$$
a = \nu^2 q^{2 z}  \; , \;\;\; b = q^{- 2 z}  \;\;\; {\rm or} \;\;\;
 a = q^{2 z} \; ,
\;\;\; b = \nu^2  q^{- 2 z}  \;
$$
for some $z \in {\mathbb Z}$.

\subsection{The case $\kappa_i V_D = 0$: Hecke algebra case \cite{IsOg}.}

{} Representations of $\alpha B\!M\!W_2$ with $\kappa_i V_D =0$ reduce to representations of the
affine Hecke algebra $\alpha H_2$. In the basis
$(e_2,e_3) = (\psi,\sigma_i \psi)$ the matrices of the generators are
 \be
 \lb{bmw12}
\sigma_i =
\left(
\begin{array}{cc}
0 & 1 \\
1 & q-q^{-1}
\end{array}
\right)  , \;
y_i =
\left(
\begin{array}{cc}
a & - (q-q^{-1}) b  \\
0 & b
\end{array}
\right) , \;
y_{i+1} =
\left(
\begin{array}{cc}
b & (q-q^{-1}) b  \\
0 & a
\end{array}
\right) ,
 \ee
where $a \neq b$ (otherwise $y_i,y_{i+1}$ are not diagonalizable). By Lemma 1, we have for
$y_i,y_{i+1} \in \alpha' B\!M\!W_2$ the eigenvalues
$a,b \in \{ q^{2 \mathbb Z}, \nu^2 q^{2 \mathbb Z} \}$. The 2-dimensional
representation (\ref{bmw12}) contains a 1-dimensional subrepresentation iff $a = q^{\pm 2} b$.
Graphically these 1- and 2-dimensional irreps of $\alpha' B\!M\!W_2$ are visualized as:

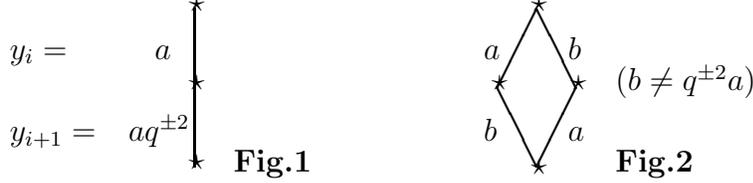
\unitlength=3.5mm
\begin{picture}(20,7.5)
\thicklines
\put(7,1){\line(0,1){6}}
\put(6.8,0.9){$\star$}
\put(6.8,6.9){$\star$}
\put(6.8,3.9){$\star$}

 \put(0,2){$y_{i+1}=$}
 \put(4.5,2){$a q^{\pm 2}$}
 \put(0,5.1){$y_i=$}
 \put(5.5,5.1){$a$}

\put(8.5,0.8){{\bf Fig.1}}

\thicklines

\put(18,2){$b$}
 \put(18,5.1){$a$}
 \put(21.2,5.1){$b$}
 \put(21.2,2){$a$}
 \put(23,3.9){$(b \neq q^{\pm2}a)$}
 \put(18.3,3.9){$\star$}
\put(21.3,3.9){$\star$} 
\put(19.8,6.9){$\star$}
\put(20,7){\line(-1,-2){1.5}} 
\put(20,7){\line(1,-2){1.5}}
\put(20,1){\line(-1,2){1.5}} 
\put(20,1){\line(1,2){1.5}}
\put(19.8,0.7){$\star$}

\put(23,0.8){{\bf Fig.2}}
\end{picture}

\vspace{-0.2cm}

\noindent
Different paths going from the upper vertex to the lower vertex correspond to different
eigenvectors of $y_i,y_{i+1}$. The indices on the edges are eigenvalues of
$y_i,y_{i+1}$.

\subsection{$\kappa_i V_D \neq 0$: odd dimensional representations for $\alpha' B\!M\!W_2$}

Using condition (\ref{bmw10}) for the reduction $V_\infty$ to $V_{2m+1}$, one can describe odd
dimensional representations of
$\alpha' B\!M\!W_2$, determine matrices for the
action of $y_i$, $y_{i+1}$ on $V_{2m+1}$ and calculate
\be
\lb{detev}
\det (y_i) = \prod_{r=1}^{2m+1} y^{(r)}_i = \nu^{2m} \; , \;\;\;
\det(y_{i+1}) = \prod_{r=1}^{2m+1} y^{(r)}_{i+1} = \nu^{2m+2} \; .
\ee
Here for eigenvalues $y^{(r)}_i$, $y^{(r)}_{i+1}$ $(r = 1,2, \dots, 2m+1)$ of
$y_i$ and $y_{i+1}$ we have constraints
$$
y^{(r)}_i \, y^{(r)}_{i+1} = \nu^2 \; , \;\;\; r=1, \dots, 2m+1 \;
$$
and (see eq. (\ref{spec2}))
 $$
y^{(r)}_i \in \{ q^{2 {\mathbb Z}} , \; \nu^2 q^{2 {\mathbb Z}} \} \; , \;\; r=1, \dots, 2m+1 \; .
  $$

These odd-dimensional irreps are visualized as graphs:

\vspace{0.2cm}

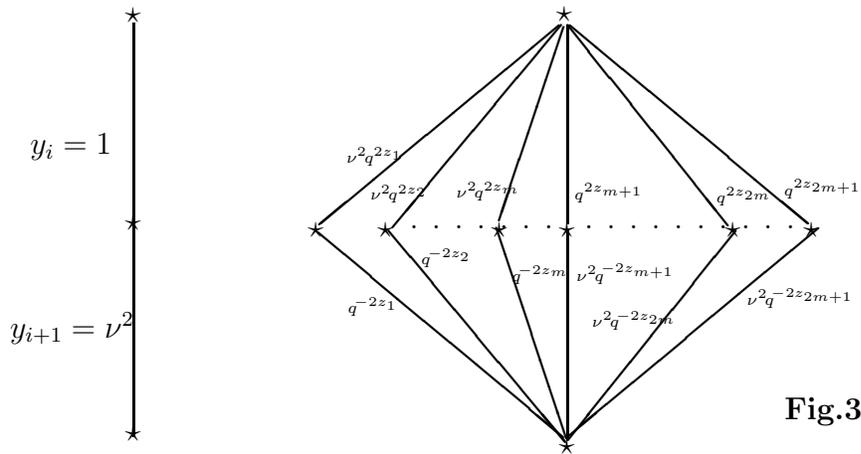
\setlength{\unitlength}{0.00016in}%

\begin{picture}(18225,15105)(1426,-15322)
%\put(29651,-8236){\tiny  $\Lambda_{+}$}
%\put(27026,-8311){\tiny  $\Lambda_{+}$}
%\put(21551,-8311){\tiny  $\Lambda_{+}$}
%\put(17751,-8236){\tiny  $\Lambda_{-}$}
%\put(14051,-8236){\tiny  $\Lambda_{-}$}
%\put(11426,-8236){\tiny $\Lambda_{-}$}
%\put(20051,-561){\tiny  $\Lambda$}
%\put(21401,-15286){\tiny $\Lambda$}

\put(14176,-6736){\tiny $\nu^2 \! q^{2z_2}$}
\put(13351,-5461){\tiny $\nu^2 \! q^{2z_1}$}
\put(17126,-6661){\tiny $\nu^2 \! q^{2z_m}$}
\put(21176,-6736){\tiny $q^{2z_{m+1}}$}
\put(26051,-6811){\tiny $q^{2z_{2m}}$}
\put(28376,-6586){\tiny $q^{2z_{2m+1}}$}
\put(27101,-10336){\tiny $\nu^2 \! q^{\! -2z_{2m+1}}$}
\put(21776,-11086){\tiny $\nu^2 \! q^{\! -2z_{2m}}$}
\put(21251,-9511){\tiny $\nu^2 \! q^{\! -2z_{m+1}}$}
\put(19026,-9511){\tiny $q^{\! -2z_{m}}$}
\put(15926,-8936){\tiny $q^{\! -2z_2}$}
\put(13401,-10611){\tiny $q^{\! -2z_1}$}

\thicklines
\put(20571,-606){$\star$}
\put(20651,-8011){$\star$}
\put(29051,-8011){$\star$}
\put(26351,-8011){$\star$}
\put(15526,-7711){. . . . . . .}
\put(21526,-7711){. . . . . . . . . .}
\put(18326,-8011){$\star$}
\put(14426,-8011){$\star$}
\put(12026,-8011){$\star$}
\put(20651,-15486){$\star$}
\put(20651,-736){\line(-6,-5){8225.410}}
\put(20726,-736){\line(-5,-6){5784.836}}
\put(20951,-736){\line( 0,-1){6825}}
\put(21026,-736){\line( 6,-5){8181.148}}
\put(20951,-736){\line( 4,-5){5495.122}}
\put(20801,-811){\line(-1,-3){2212.500}}
\put(12376,-7861){\line( 6,-5){8424.590}}
\put(14801,-7786){\line( 5,-6){6061.475}}
\put(18501,-7711){\line( 1,-3){2407.500}}
\put(20951,-7786){\line( 0,-1){7200}}
\put(26726,-7786){\line(-4,-5){5765.854}}
\put(29551,-7711){\line(-6,-5){8535.246}}

\put(5751,-656){$\star$}
\put(5751,-7786){$\star$}
\put(5751,-15036){$\star$}
\put(6051,-656){\line( 0,-1){14200}}

\put(2501,-5161){$y_i=1$}
\put(1801,-11461){$y_{i+1}=\nu^2$}

\put(28401,-14286){\bf Fig.3}

\end{picture}

\noindent
where $z_r \in {\mathbb Z}$ and $\sum_{r=1}^{2m+1} z_r =0$ as it follows from (\ref{detev}). 
Different paths going from the top vertex to the bottom vertex correspond to different
 common eigenvectors of 
$y_i,y_{i+1}$. Indices on upper and lower edges of these paths are the eigenvalues of $y_i$ and
$y_{i+1}$, respectively.

\noindent
{\bf Remark.} In view of the braid relations $\sigma_{i} \sigma_{i \pm 1} \sigma_i =
\sigma_{i \pm 1} \sigma_i  \sigma_{i \pm 1}$ and
possible eigenvalues of $\sigma$'s for 1-dimensional representations  (described in Subsections 3.2
and 3.3) we conclude that the following chains of 1-dimensional representations are forbidden

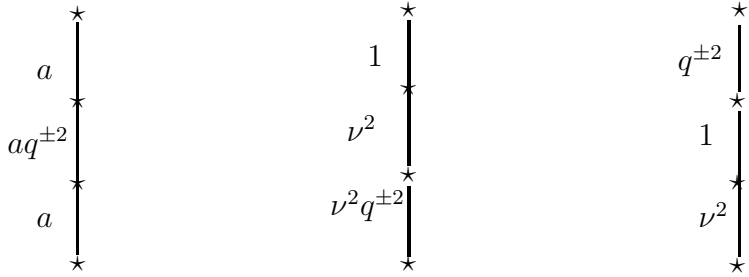
\unitlength=5.5mm
\begin{picture}(20,7.5)
\thicklines
\put(5,1.15){\line(0,1){5.6}}

\put(4.8,0.75){$\star$}
\put(4.8,6.8){$\star$}
\put(4.8,2.7){$\star$}
\put(4.8,4.7){$\star$}

\put(4,1.8){$a$}
\put(3.3,3.6){$a q^{\pm 2}$}
\put(4,5.4){$a$}

\put(13,1.1){\line(0,1){1.7}}
\put(13,3.3){\line(0,1){3.5}}
\put(12.8,0.75){$\star$}
\put(12.8,6.9){$\star$}
\put(12.8,2.9){$\star$}
\put(12.8,5){$\star$}
\put(11.1,2.1){$\nu^2 q^{\pm 2}$}
\put(11.5,3.9){$\nu^{2}$}
\put(12,5.7){$1$}

\put(21,1.1){\line(0,1){3.5}}
\put(21,5.1){\line(0,1){1.6}}
\put(20.75,0.7){$\star$}
\put(20.8,6.9){$\star$}
\put(20.75,2.7){$\star$}
\put(20.75,4.7){$\star$}
\put(20,1.8){$\nu^2 $}
\put(20,3.8){$1$}
\put(19.5,5.6){$q^{\pm 2}$}

\end{picture}

\vspace{0.5cm}
\noindent
where $a= q^{2z}$ or $a = \nu^2 q^{2z}$ $(z \in {\mathbb Z})$.

\subsection{$\kappa_i V_D \neq 0$: even dimensional representations of $\alpha' B\!M\!W_2$}

With the help of conditions (\ref{bmw11}) we reduce $V_\infty$ to $V_{2m}$, then explicitly
construct
$(2m) \times (2m)$ matrices for the operators $y_i,y_{i+1}$ and calculate their determinants
\be
\lb{bmw14}
\det(y_i) =
\prod_{r=1}^{2m} y_i^{(r)} =  \epsilon q^{\epsilon} \, \nu^{2m-1}  \; , \;\;
\det(y_{i+1}) =
\prod_{r=1}^{2m} y_{i+1}^{(r)} = - \epsilon q^{\epsilon} \, \nu^{2m+1}
\; ,
\ee
where $y_i^{(r)}$, $y_{i+1}^{(r)}$ are eigenvalues of $y_i$, $y_{i+1}$ (we have two possibilities: 
$\epsilon=\pm 1$). We see from
(\ref{bmw14}) that all
$(2m)$ eigenvalues of $y_i$, $y_{i+1}$ cannot belong to the spectrum  (\ref{spec2}).
More precisely there is at least one eigenvalue $y_i^{(r)}$ of $y_i$
(and the eigenvalue $y_{i+1}^{(r)}$ of $y_{i+1}$) such that
$$
y_i^{(r)},y_{i+1}^{(r)} \notin \{ q^{2 {\mathbb Z}} , \; \nu^2 q^{2 {\mathbb Z}} \} \; .
$$
Thus, even dimensional irreps of
$\alpha B\!M\!W_2$ subject to the conditions (\ref{bmw11}),
are not admissible for $\alpha' B\!M\!W_2$.

\section{\bf Representations of $B\!M\!W$ algebras}

\subsection{Spec$(y_1,\dots,y_{n})$ and rules for strings of eigenvalues}

Now we reconstruct the representation theory of $B\!M\!W$ algebras using an approach which
generalizes the approach of Okounkov -- Vershik \cite{OV} for symmetric groups.

The  JM elements {$\{
\tilde{y}_1,\dots,\tilde{y}_{n}\}$} generate a commutative
subalgebra in $B\!M\!W_n$. The basis in the space of an irrep of
$B\!M\!W_n$ can be chosen to be the common eigenbasis of all
$\tilde{y}_i$. Each common  eigenvector $v$ of $\tilde{y}_i$,
$$
\tilde{y}_i \,  v = a_i \, v \; , \;\; i=1,\dots,n \; ,
$$
defines a string $(a_1, \dots,a_n) \in {\mathbb C}^n$. Denote by ${\rm
Spec}(\tilde{y}_1,\dots,\tilde{y}_{n})$ the set of such strings.

We summarize our results about representations of $\alpha' B\!M\!W_2$ and the spectrum of the JM
elements $\tilde{y}_i$ in the following Proposition.

\vspace{0.2cm}

\noindent{\bf Proposition 4.} {\it Consider the string
$$
\alpha = (a_1,...,a_i,a_{i+1},...,a_n) \in  {\rm
Spec}(\tilde{y}_1,...,\tilde{y}_i,\tilde{y}_{i+1},...,\tilde{y}_{n}) \; .
$$
Let $v_\alpha$ be the corresponding eigenvector of $\tilde{y}_i$:
$\tilde{y}_i \, v_\alpha = a_i \, v_\alpha$. Then
$$
\!\!\!\!\!\!\! \begin{array}{l}
 {\bf (1)} \;\;\;\; a_i \in \{ q^{2 {\mathbb Z}} ,
\; \nu^2 q^{2 {\mathbb Z}} \} \, ; \\
 {\bf (2)} \;\;\;\; a_i \neq a_{i+1} \; , \;\; i =1,\dots,n-1 \; ; \\
 {\bf (3a)} \;\; a_i a_{i+1} \neq \nu^2 \; ,
 \;\;\; a_{i+1} = q^{\pm 2} \, a_{i} \Rightarrow
\sigma_i \cdot v_\alpha = \pm q^{\pm 1} v_\alpha \, , \;\;
\kappa_i \cdot v_\alpha = 0 \, ; \\ [0.2cm]
 {\bf (3b)} \;\; a_i a_{i+1} \neq \nu^2 \, ,
  \;\;\; a_{i+1} \neq q^{\pm 2} \, a_{i} \Rightarrow \\
\;\;\;\;\;\;\;\;\alpha\, '= (a_1,..., a_{i+1},a_i,...,a_n) \in  {\rm
Spec}(\tilde{y}_1,...,\tilde{y}_i,\tilde{y}_{i+1},...,\tilde{y}_{n})
\, , \; \kappa_i \cdot v_\alpha = 0 \, , \; \kappa_i \cdot v_{\alpha\, '} = 0 \, ; \\ [0.25cm]
 {\bf (4)} \;\;\;\; a_i
a_{i+1} \! = \! \nu^2 \;\; \Rightarrow \;\; \exists \;\; {\rm odd} \; {\rm number} \; {\rm of} \;
{\rm strings} \;\;
\alpha^{(k)} \;\; (k=1,2,\dots, 2m+1): \\ [0.2cm]
\;\;\;\;\;\;\;\;\;
 \alpha^{(k)}=(a_1,...,a_{i-1},a^{(k)}_{i}\!\!
,a^{(k)}_{i+1},a_{i+2},...,a_n) \in {\rm Spec}(\tilde{y}_1,...,\tilde{y}_n) \;\; \forall k \; ,
\\ [0.3cm]
\;\;\;\;\;\;\;\;
\;\; \alpha \in \{\alpha^{(k)}\} \, , \;\; a^{(k)}_{i} a^{(k)}_{i+1}
\! = \! \nu^2 \, , \;\; \prod\limits_{k=1}^{2m+1} a^{(k)}_{i} =
\nu^{2m} \, , \;\; \prod\limits_{k=1}^{2m+1} a^{(k)}_{i+1} =
\nu^{2m+2} \; .
\end{array}
$$}

The necessary and sufficient conditions for a string to belong to the common spectrum of
$\tilde{y}_i$ are formulated in the following way.

\vspace{0.2cm}

\noindent
{\bf Proposition 5.} {\it The string $\alpha =(a_1,a_2, \dots , a_n)$, where $a_i \in
(q^{2\mathbb Z}, \nu^2 q^{2\mathbb Z})$, belongs to the set
${\rm Spec}(\tilde{y}_1,\tilde{y}_2,...,\tilde{y}_{n})$
 iff $\alpha$ satisfies the following conditions ($z \in {\mathbb Z}$)
$$
\begin{array}{l}
{\bf (1)} \;\; a_1 = 1  \; ; \;\;\;\;\;\; {\bf (2)} \;\;\;\; a_i = \nu^2 q^{-2 z}
\;\;
\Rightarrow q^{2z} \in \{a_1, ... , a_{i-1} \} \; ; \\ [0.2cm]
{\bf (3)} \;\; a_i = q^{2z} \Rightarrow
\{ a_i q^2, a_i q^{-2} \} \cap \{a_1, ... , a_{i-1} \} \neq \O
 \; , \;\;  z \neq 0 ; \\ [0.2cm]
\!\!\!\!\!\! {\bf (4a)} \;\; a_i \! = \! a_j \! = \! q^{2z} \; (i < j) \! \Rightarrow
\!\! \left\{ \!
\begin{array}{l}
\!\! {\rm either} \; \{ q^{^{2(z + 1)}},q^{^{2(z - 1)}} \} \! \subset \!
\{a_{_{i+1}}, ... , a_{_{j-1}} \} \; , \\ [0.2cm]
\!\! {\rm or} \; \nu^2 q^{-2z}
\in \{a_{_{i+1}}, ... , a_{_{j-1}} \} \; ;
\end{array}
\right.
\\ [0.6cm]
\!\!\!\!\!\! {\bf (4b)} \; a_i \! = \! a_j \! = \! \nu^2 q^{2z} \, (i < j) \Rightarrow
\!\! \left\{ \!
\begin{array}{l}
\!\! {\rm either} \; \{ \nu^2 q^{^{2(z + 1)}} \!\! , \nu^2 q^{^{2(z - 1)}} \} \! \subset \!
\{a_{_{i+1}}, ... , a_{_{j-1}} \} \; , \\ [0.2cm]
\!\! {\rm or} \;  q^{-2z}
\in \{a_{_{i+1}}, ... , a_{_{j-1}} \} \; ;
\end{array}
\right.
\\ [0.5cm]
\!\!\!\!\!\! {\bf (5a)} \;\; a_i = \nu^2 q^{-2z} \; , \;\; a_j = q^{2 z'} \; (i < j)
\Rightarrow
q^{2z} \; {\rm or} \; \nu^2 q^{-2 z'}\in \{a_{_{i+1}}, ... , a_{_{j-1}} \}  \; ; \\ [0.5cm]
\!\!\!\!\!\! {\bf (5b)} \;\; a_i = q^{2z} \; , \;\; a_j = \nu^2 q^{-2 z'} \; (i < j)
\Rightarrow
\nu^2 q^{-2z} \; {\rm or} \;  q^{2 z'} \in \{a_{_{i+1}}, ... , a_{_{j-1}} \}  \; .
\end{array}
$$
where in {\bf (5a)} and {\bf (5b)} we set $z' = z \pm 1$}.

\subsection{Young graph for $B\!M\!W$ algebras}

We illustrate the above considerations on the example of the colored (in the sense of
\cite{IsOg}) Young graph for the algebra
$BMW_5$. This graph contains the whole information about the irreps of $BMW_5$ and the branching rules
$BMW_5 \downarrow BMW_4$.

\setlength{\unitlength}{0.00023495in}%
\begingroup\makeatletter\ifx\SetFigFont\undefined%

\fi\endgroup%
\begin{picture}(15266,24019)(4668,-23168)

 \put(23801,339){$\emptyset$}
\put(24151,-436){\tiny $1$}
\put(23476,-3061){\circle*{300}}
\put(23101,-3061){\circle*{300}}
 \put(22951,-2236){\tiny $q^2$}
  \put(22126,-5311){\tiny $q^4$}
    \put(21901,-10636){\tiny $q^{\!-2}$}
 \put(20626,-20361){\tiny $1$}
 \put(22501,-6361){\circle*{300}}
 \put(22126,-6361){\circle*{300}} \put(22876,-6361){\circle*{300}}
\put(24001,-1111){\circle*{300}}

\put(21001,-22186){\circle*{300}} \put(20626,-22186){\circle*{300}}
\put(20251,-22186){\circle*{300}} \put(20251,-22561){\circle*{300}}
\put(20626,-22561){\circle*{300}}

 \put(22501,-11836){\circle*{300}}
 \put(22501,-11461){\circle*{300}}
  \put(22876,-11461){\circle*{300}}
\put(23251,-11461){\circle*{300}}

 \put(22501,-6661){\line(0,-1){4500}}
  \put(24001,239){\line( 0,-1){1050}}
 \put(23251,-2761){\line( 2, 5){600}}

 \put(23176,-3436){\line(-1,-4){675}}

 \put(23026,-11986){\line(-1,-4){2452.941}}

  \put(17026,-7836){\tiny $q^{-2}$}
 \put(15301,-9361){\tiny $q^2$}

 \put(21226,-10261){\tiny $q^6$}
 \put(25876,-4636){\tiny $q^{-4}$}
 \put(24901,-4936){\tiny $q^2$}
   \put(23476,-5911){\tiny $q^{-2}$}

 \put(27301,-8386){\tiny $q^{-6}$}
 \put(24526,-9961){\tiny $1$}
 \put(25201,-8761){\tiny $q^{-4}$}
 \put(26476,-9361){\tiny $q^2$}
 \put(23476,-9436){\tiny $q^4$}

\put(24901,-2236){\tiny $q^{-2}$}
\put(16751,-4561){\tiny $1$}
\put(18301,-20461){\tiny $q^{\!-2}$}
\put(17026,-21136){\tiny $q^8$}
\put(15101,-20861){\tiny $q^{\!-4}$}
\put(13176,-16861){\tiny $q^{\!-2}$}
\put(23551,-20986){\tiny $q^4$}
\put(28076,-16411){\tiny $q^{\!-8}$}
\put(27201,-17461){\tiny $q^2$}
\put(26506,-14886){\tiny $q^{\!-6}$}

\put(21451,-20536){\tiny $q^4$}
\put(22676,-18586){\tiny $q^{\!-4}$}
\put(19651,-20461){\tiny $q^6$}
\put(25801,-17011){\tiny $1$}
\put(24651,-18151){\tiny $q^{\!\!-4}$}
\put(15201,-14986){\tiny $q^2$}
\put(10106,-20556){\tiny $q^{4}$}
\put(10201,-16036){\tiny $1$}

 \thicklines
 \put(24451,-12286){\line(-2,-5){3791.379}}
   \put(17251,-11536){\line(-1,-2){4860}}
 \put(15151,-11311){\line(-1,-2){5085}}

 \put(26176,-12736){\line(-1,-3){3000}}
  \put(24901,-3586){\line(-1,-6){450}}

 \put(24451,-6961){\line(2,-5){1800}}
 \put(26326,-6961){\line( 0,-1){4425}}
 \put(24301,-6961){\line( 0,-1){4425}}
 \put(24151,-6961){\line(-1,-3){1350}}

 \put(22276,-6586){\line(-1,-3){1425}}
 \put(16576,-6586){\line( 1,-6){675}}

 \put(16351,-6586){\line(-1,-4){1050}}
 \put(25051,-3511){\line( 1,-2){1200}}

 \put(23326,-3286){\line(1,-3){975}}
   \put(26326,-12661){\line(-1,-6){1485.811}}

 \put(24601,-12286){\line(0,-1){9150}}
 \put(26401,-12361){\line( 0,-1){9000}}
 \put(27901,-12136){\line(-1,-6){1526.351}}
 \put(27976,-12136){\line(0,-1){9075}}
 \put(16501,-3361){\line( 0,-1){2700}}

   \put(24151,-1261){\line( 2,-5){600}}

 \put(23176,-11911){\line( 0,-1){9750}}

 \put(22801,-11911){\line(-2,-5){3987.931}}
 \put(20626,-11461){\line(-1,-5){2082.692}}
 \put(20176,-11461){\line(-1,-3){3442.500}}
 \put(17326,-11611){\line(-1,-4){2497.059}}

 \put(15451,-11461){\line(-1,-3){3292.500}}
 \put(12451,-11461){\line(-2,-5){3946.552}}
 \put(26626,-6811){\line(1,-3){1200}}

 \put(12401,-11301){$\emptyset$}
 \put(16301,-3161){$\emptyset$}

  \put(21226,-11236){\circle*{300}}
 \put(24301,-11686){\circle*{300}}
\put(24301,-12061){\circle*{300}} \put(24676,-11686){\circle*{300}}
 \put(24676,-12061){\circle*{300}}
 \put(26101,-11761){\circle*{300}}
\put(26101,-12136){\circle*{300}} \put(26101,-12511){\circle*{300}}
\put(26476,-11761){\circle*{300}} \put(27901,-10636){\circle*{300}}
\put(27901,-11011){\circle*{300}} \put(27901,-11386){\circle*{300}}
\put(27901,-11761){\circle*{300}} \put(9826,-21586){\circle*{300}}
\put(12676,-21661){\circle*{300}} \put(28276,-21511){\circle*{300}}
\put(28276,-21886){\circle*{300}} \put(20851,-11236){\circle*{300}}

   \put(24901,-2911){\circle*{300}}
\put(24901,-3286){\circle*{300}}

\put(24301,-6361){\circle*{300}} \put(24676,-6361){\circle*{300}}
\put(24301,-6736){\circle*{300}} \put(26401,-6361){\circle*{300}}
\put(26401,-5986){\circle*{300}} \put(26401,-6736){\circle*{300}}
\put(16501,-6361){\circle*{300}} \put(15751,-11161){\circle*{300}}
\put(15376,-11161){\circle*{300}} \put(17401,-11011){\circle*{300}}
\put(17401,-11386){\circle*{300}} \put(20101,-11236){\circle*{300}}
\put(20476,-11236){\circle*{300}} \put(28276,-22261){\circle*{300}}
\put(19126,-22186){\circle*{300}} \put(18001,-22186){\circle*{300}}
\put(18001,-22561){\circle*{300}} \put(17026,-22036){\circle*{300}}
\put(16651,-22036){\circle*{300}} \put(16276,-22036){\circle*{300}}
\put(15901,-22036){\circle*{300}} \put(15526,-22036){\circle*{300}}
\put(14476,-21811){\circle*{300}} \put(14476,-22186){\circle*{300}}
\put(14476,-22561){\circle*{300}} \put(12301,-21661){\circle*{300}}
\put(12301,-22036){\circle*{300}} \put(10201,-21586){\circle*{300}}
\put(10576,-21586){\circle*{300}} \put(8326,-21586){\circle*{300}}
\put(18376,-22186){\circle*{300}} \put(28276,-22636){\circle*{300}}
\put(28276,-23011){\circle*{300}} \put(26326,-21661){\circle*{300}}
\put(24601,-21736){\circle*{300}} \put(24601,-22111){\circle*{300}}
\put(24226,-21736){\circle*{300}} \put(24226,-22111){\circle*{300}}
\put(24226,-22486){\circle*{300}} \put(23101,-21961){\circle*{300}}
 \put(22726,-21961){\circle*{300}} \put(22351,-21961){\circle*{300}}
 \put(22351,-22336){\circle*{300}} \put(22351,-22711){\circle*{300}}

 \put(26326,-22036){\circle*{300}}
\put(26326,-22411){\circle*{300}} \put(26326,-22786){\circle*{300}}
\put(26701,-21661){\circle*{300}} \put(18751,-22186){\circle*{300}}

 \put(24676,-3586){\line(-3,-1){7852.500}}
  \put(24076,-6811){\line(-5,-3){6529.412}}
   \put(26326,-12661){\line(-5,-4){11485.811}}
 \put(24076,-12136){\line(-6,-5){11161.475}}
 \put(22276,-11836){\line(-1,-1){9525}}
  \put(26032,-12820){\line(-3,-2){13309.615}}
 \put(20251,-11461){\line(-1,-1){9780}}
  \put(17776,-10861){\line( 2,1){8400}}
 \put(12676,-10861){\line( 5, 6){3750}}
 \put(15676,-10861){\line( 2, 1){8400}}
 \put(15451,-10861){\line( 3,2){6525}}
  \put(16651,-5986){\line( 5, 2){6375}}
   \put(22201,-11836){\line(-6,-5){11419.672}}
    \put(16976,-11461){\line(-5,-6){8250}}
 \put(15076,-11386){\line(-2,-3){6600}}
  \put(27976,-12136){\line( -4,-3){13075}}
  \put(16801,-2911){\line( 4,1){6900}}

\put(19051,-4411){\tiny $\nu^2q^{-2}$}
\put(15726,-21536){\tiny $\nu^2\!q^{6}$}
\put(16001,-20361){\tiny $\nu^2\!q^{\!-2}$}
\put(13701,-19636){\tiny $\nu^2\!q^{\!\!-4}$}
 \put(18976,-10536){\tiny $\nu^2 q^4$}
 \put(17576,-10356){\tiny $\nu^2\! q^{\!-2}$}
 \put(20246,-8236){\tiny $\nu^2q^2$}
 \put(18376,-8236){\tiny $\nu^2q^{\!-4}$}
 \put(14101,-8311){\tiny $\nu^2$}
 \put(19651,-5611){\tiny $\nu^2q^2$}
 \put(19876,-1786){\tiny $\nu^2$}
 \put(13351,-21486){\tiny $\nu^2q^4$}
\put(11326,-21136){\tiny $\nu^2q^2$}
\put(11651,-19300){\tiny $\nu^2\!q^{\!-6}$}
\put(9026,-21211){\tiny $\nu^2\!\!q^2$}
\put(11401,-15661){\tiny $\nu^2\!q^{\!-2}$}
\put(16576,-18661){\tiny $\nu^2$}

\put(4551,-20786){{\bf Fig.4}}
\end{picture}

\noindent
A vertex $\{\lambda;5\}$ on the lowest level of this graph is labeled by some Young diagram
$\lambda$; this vertex corresponds to the irrep $W_{\{\lambda;5\}}$ of
$BMW_5$ (the notation $\{\lambda;5\}$ is designed to encode the diagram $\lambda$ and the level
on which this diagram is located; the levels are counted starting from 0). Paths going down from
the top vertex $\emptyset$ to the lowest level
(that is, paths of length 5) correspond to common eigenvectors of the JM elements $\tilde{y}_1,
\dots , \tilde{y}_5$. Paths ending at $\{\lambda;5\}$ label the basis in $W_{\{\lambda;5\}}$.
In particular, the number of different paths going down from the top $\emptyset$ to
$\{\lambda;5\}$ is equal to the dimension
of the irrep $W_{\{\lambda;5\}}$.

Note that the colored Young graph in Fig.4 contains subgraphs presented in Fig.1, Fig.2 and Fig.3. 
For example, in Fig.4 one recognizes rhombic subgraphs (the vertices on the subgraphs are obtained 
from one another by a rotation) 

\vspace{0.3cm}

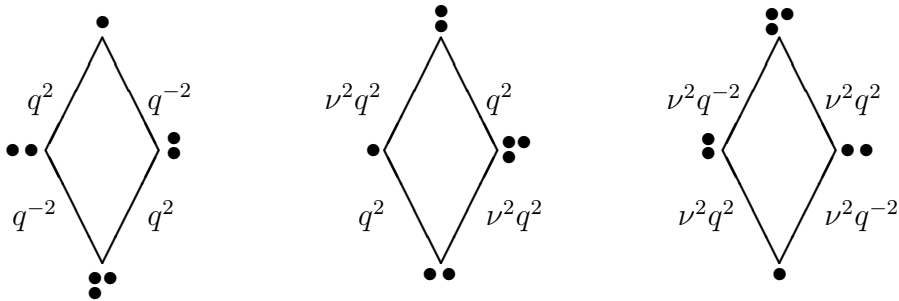
\unitlength=5mm
\begin{picture}(20,7.5)

\thicklines

\put(1.6,2){$q^{-2}$}
 \put(2,5.1){$q^{2}$}
 \put(5.2,5.1){$q^{-2}$}
 \put(5.2,2){$q^{2}$}

\put(5.7,4.1){$\bullet$}
\put(5.7,3.7){$\bullet$}

\put(1.4,3.8){$\bullet$}
\put(1.9,3.8){$\bullet$}

\put(3.6,0.4){$\bullet$}
\put(3.6,0.0){$\bullet$}
\put(4,0.4){$\bullet$}

\put(4,7){\line(-1,-2){1.5}} 
\put(4,7){\line(1,-2){1.5}}
\put(4,1){\line(-1,2){1.5}} 
\put(4,1){\line(1,2){1.5}}
\put(3.8,7.2){$\bullet$}

%%%%%%%%%%%%%%%%%%%%%%%%%%%%%%%

\put(10.8,2){$q^{2}$}
 \put(9.9,5.1){$\nu^2 q^{2}$}
 \put(14.2,5.1){$q^{2}$}
 \put(14.2,2){$\nu^2 q^{2}$}
 \put(11,3.8){$\bullet$}
 
\put(14.6,4){$\bullet$}
\put(14.6,3.6){$\bullet$}
\put(15,4){$\bullet$} 

\put(12.5,0.5){$\bullet$}
\put(13,0.5){$\bullet$}

\put(13,7){\line(-1,-2){1.5}} 
\put(13,7){\line(1,-2){1.5}}
\put(13,1){\line(-1,2){1.5}} 
\put(13,1){\line(1,2){1.5}}

\put(12.8,7.5){$\bullet$}
\put(12.8,7.1){$\bullet$}

%%%%%%%%%%%%%%%%%%%%%%%%%%%

\put(19.3,2){$\nu^2 q^{2}$}
 \put(19,5.1){$\nu^2 q^{-2}$}
 \put(23.2,5.1){$\nu^2 q^{2}$}
 \put(23.2,2){$\nu^2 q^{-2}$}

\put(19.9,4.1){$\bullet$}
\put(19.9,3.7){$\bullet$}

\put(23.6,3.8){$\bullet$}
\put(24.1,3.8){$\bullet$}

\put(21.6,7.4){$\bullet$}
\put(21.6,7.0){$\bullet$}
\put(22,7.4){$\bullet$}

\put(22,7){\line(-1,-2){1.5}} 
\put(22,7){\line(1,-2){1.5}}
\put(22,1){\line(-1,2){1.5}} 
\put(22,1){\line(1,2){1.5}}
\put(21.8,0.5){$\bullet$}

\end{picture} 

\vspace{0.3cm}

\noindent
of the type presented in Fig.2.

Let $(s,t)$ be coordinates of a node in the Young diagram $\lambda$. To the node $(s,t)$ of the
diagram $\lambda$ we associate a number
$q^{2(s-t)}$ which is called "content":

\vspace{0.5cm}

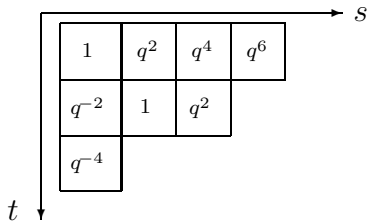
\unitlength=5mm
\begin{picture}(25,4)
\put(7,4){\vector(1,0){8}}
\put(7,4){\vector(0,-1){5.5}}
\put(15.3,3.8){$s$}
\put(6.1,-1.5){$t$}

\put(7.5,1.3){$
\begin{tabular}{|c|c|c|c|}
\hline
  $\!\!\! $  $_1$ & $\!\! $  $_{q^2}$  &
  $\!\! $  $_{q^4}$ & $\!\! $  $_{q^6}$ \\[0.2cm]
\hline
  $\!\!\! $ $ _{q^{\! -2}}$  &   $\!\!\!$ $_1$  & $\!\!\! $ $_{q^2}$  &
  \multicolumn{1}{c}{} \\[0.2cm]
\cline{1-3}
  $\!\!\! $ $ _{q^{\! -4}}$  & \multicolumn{1}{c}{} \\[0.2cm]
\cline{1-1}
\end{tabular}
$}
\end{picture}

\vspace{1cm}

\noindent
Then according to the colored Young graph on Fig.4, at each step down along the path one can add or 
remove one node (therefore this graph is called the "oscillating" Young graph) and the eigenvalue 
of the corresponding JM element is determined by the content of the node:

\vspace{0.2cm}

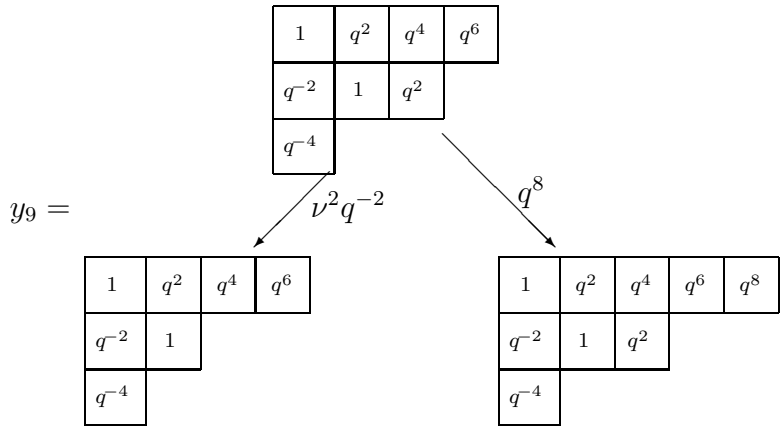
\unitlength=5mm
\begin{picture}(25,4)

\put(7.5,1.3){$
\begin{tabular}{|c|c|c|c|}
\hline
  $\!\!\! $  $_1$ & $\!\! $  $_{q^2}$  &
  $\!\! $  $_{q^4}$ & $\!\! $  $_{q^6}$ \\[0.2cm]
\hline
  $\!\!\! $ $ _{q^{\! -2}}$  &   $\!\!\!$ $_1$  & $\!\!\! $ $_{q^2}$  &
  \multicolumn{1}{c}{} \\[0.2cm]
\cline{1-3}
  $\!\!\! $ $ _{q^{\! -4}}$  & \multicolumn{1}{c}{} \\[0.2cm]
\cline{1-1}
\end{tabular}
$}
\end{picture}

\vspace{1.3cm}

\unitlength=5mm
\begin{picture}(25,4)

\put(9,6){\vector(-1,-1){2}}
\put(12,7){\vector(1,-1){3}}

\put(0.5,4.8){$y_9 =$}
\put(8.5,4.8){$\nu^2 q^{-2}$}
\put(14,5.2){$q^8$}

\put(2.5,1.3){$
\begin{tabular}{|c|c|c|c|}
\hline
  $\!\!\! $  $_1$ & $\!\! $  $_{q^2}$  &
  $\!\! $  $_{q^4}$ & $\!\! $  $_{q^6}$ \\[0.2cm]
\hline
  $\!\!\! $ $ _{q^{\! -2}}$  &   $\!\!\!$ $_1$  &
  \multicolumn{1}{c}{} \\[0.2cm]
\cline{1-2}
  $\!\!\! $ $ _{q^{\! -4}}$  & \multicolumn{1}{c}{} \\[0.2cm]
\cline{1-1}
\end{tabular}
$}

\put(13.5,1.3){$
\begin{tabular}{|c|c|c|c|c|}
\hline
  $\!\!\! $  $_1$ & $\!\! $  $_{q^2}$  &
  $\!\! $  $_{q^4}$ & $\!\! $  $_{q^6}$ & $\!\! $  $_{q^8}$ \\[0.2cm]
\hline
  $\!\!\! $ $ _{q^{\! -2}}$  &   $\!\!\!$ $_1$  & $\!\!\! $ $_{q^2}$  &
  \multicolumn{1}{c}{} \\[0.2cm]
\cline{1-3}
  $\!\!\! $ $ _{q^{\! -4}}$  & \multicolumn{1}{c}{} \\[0.2cm]
\cline{1-1}
\end{tabular}
$}
\end{picture}

\vspace{1cm}

\noindent
The eigenvalue corresponding to the addition or removal of the $(s,t)$ node is $q^{2(s-t)}$ or
$\nu^2 q^{-2(s-t)}$, respectively.

 Let $X(n)$ be the set of paths of length $n$ starting from the top vertex $\emptyset$ and going down in
the Young graph of oscillating Young diagrams. Now we formulate the following Proposition.

\vspace{0.1cm}

\noindent
 {\bf Proposition 6.} {\it
There is a bijection between the set $\, {\mathrm Spec} (\tilde{y}_1, \dots , \tilde{y}_n)$ and the
set
$X(n)$.}

\subsection{Primitive idempotents}

The colored Young graph (as in Fig. 4) gives also the rule of construction of a complete set of
 orthogonal primitive idempotents for the
$BMW$ algebra. The completeness of the set of orthogonal primitive idempotents is equivalent to the
maximality of the commutative set of JM elements. Let $\{\lambda;n\}$ be a vertex in the Young
graph with

{\unitlength=4mm
\begin{picture}(25,6.5)
\put(0,3){$\lambda = $} \put(4,5.5){\line(1,0){5}}
\put(4,4){\line(1,0){5}} \put(4,4){\line(0,1){1.5}}
\put(9,4){\line(0,1){1.5}}

\put(4,2.5){\line(1,0){3}} \put(7,2.5){\line(0,1){1.5}}
\put(4,4){\line(1,0){3}} \put(4,1){\line(0,1){3}}
\put(6,1.5){\line(0,1){1}} \put(4,1.5){\line(1,0){2}}
\put(4.5,1){$\dots$}

\put(4,0){\line(0,1){0.5}} \put(5,0){\line(0,1){0.5}}
\put(4,0){\line(1,0){1}} \put(4,0.5){\line(1,0){1}}

\put(6,6.2){$_{_{\lambda_{_{(1)}}}}$} \put(3,4.7){$_{n_{_1}}$}
\put(9.1,4){$_{n_{_1},} {_{\lambda_{_{(1)}}}}$}
\put(7.1,2.5){$_{n_{_2},} {_{\lambda_{_{(2)}}}}$}
\put(6.2,1.5){$_{n_{_3},} {_{\lambda_{_{(3)}}}}$}
\put(5.1,0){$_{n_{_k},} {_{\lambda_{_{(k)}}}}$}

 \put(13,3.4){$(n_i,\lambda_{(i)})$ are coordinates of the nodes}
\put(13.5,2.2){which are in the corners of
$\lambda=[\lambda_{(1)}^{n_1},\lambda_{(2)}^{n_2-n_1},\dots,\lambda_{(k)}^{n_k-n_{k-1}}]$.}

\end{picture}}

\vspace{0.1cm}

\noindent Consider any path  $T_{\{\lambda;n\}}$ going down from the top $\emptyset$ to this vertex. Let
$E_{T_{\{\lambda;n\}}} \in B\! M\! W_{n}$ be the primitive
idempotent corresponding to $T_{\{\lambda;n\}}$. Using the branching rule implied by the Young
graph for
 $B\! M\! W_{n+1}$ we know all possible eigenvalues of the
 element $\tilde{y}_{n+1}$ and, therefore, obtain the identity
$$
E_{T_{\{\lambda;n\}}} \cdot \prod_{r=1}^{k+1} \left(
\tilde{y}_{n+1} - q^{2(\lambda_{(r)} - n_{r-1})} \right)
\prod_{r=1}^{k} \left( \tilde{y}_{n+1} - \nu^2 q^{2(n_{r} -
\lambda_{(r)})} \right)= 0 \; ,
$$
where  $\lambda_{(k+1)} = n_0 = 0$. So, for a new diagram $\lambda'$ obtained by adding to
$\lambda$ a new node with coordinates $(n_{j-1} +1,
\lambda_{(j)}+1)$ the corresponding primitive idempotent (after an appropriate normalization) reads
$$
E_{T_{\{\lambda';n+1\}}} =  E_{T_{\{\lambda;n\}}} \cdot \prod_{\stackrel{r=1}{_{r \neq j}}}^{k+1}
\! \frac{\left( \tilde{y}_{_{n+1}} - q^{2(\lambda_{(r)} -
n_{r-1})} \right)}{\left( q^{2(\lambda_{(j)} - n_{j-1})} -
q^{2(\lambda_{(r)} - n_{r-1})} \right)}
 \prod_{r=1}^{k} \! \frac{\left(
\tilde{y}_{_{n+1}} - \nu^2 q^{2(n_{r} - \lambda_{(r)} )}
\right)}{\left( q^{2(\lambda_{(j)} - n_{j-1})} - \nu^2 q^{2(n_{r} -
\lambda_{(r)} )} \right)}\,  .
 $$
For a new diagram $\lambda''$ which is obtained from
$\lambda$ by removing a node with coordinates
$(n_{j},\lambda_{(j)})$ we construct the primitive idempotent
$$
E_{T_{\{\lambda'';n+1\}}} =  E_{T_{\{\lambda;n\}}} \cdot \prod_{r=1}^{k+1}
\! \frac{\left( \tilde{y}_{_{n+1}} - q^{2(\lambda_{(r)} -
n_{r-1})} \right)}{\left( \nu^2 q^{2(n_{j}-\lambda_{(j)})} - q^{2(\lambda_{(r)} - n_{r-1})}
\right)}
 \prod_{\stackrel{r=1}{_{r \neq j}}}^{k} \! \frac{\left(
\tilde{y}_{_{n+1}} - \nu^2 q^{2(n_{r} - \lambda_{(r)} )}
\right)}{\left( \nu^2 q^{2(n_{j}-\lambda_{(j)})} - \nu^2 q^{2(n_{r} -
\lambda_{(r)} )} \right)}\,  .
 $$
 Using these
formulas and the "initial data" $E_{T_{\{\emptyset;0\}}} =1$, one can deduce step by step explicit
expressions for the primitive orthogonal idempotents related to the paths in the BMW Young graph.

\section{\bf Outlook}

In this paper we reconstructed the representation theory of the tower of the $BMW$ algebras using 
the properties of the commutative subalgebras, generated by the Jucys--Murphy elements, in the 
$BMW$ algebras. This representation theory is of use in the representation theory of the quantum 
groups 
$U_q(osp(N|K))$ due to the Brauer - Schur - Weyl duality, but
 also find applications in physical models.
Recently \cite{IsOg3} we have formulated integrable chain models with nontrivial boundary 
conditions in terms of the affine Hecke algebras $\alpha H_n$ and the affine BMW algebras $\alpha 
BMW_n$. The Hamiltonians for these models are special elements of the algebras $\alpha H_n$ and 
$\alpha BMW_n$. E.g., for the $\alpha BMW_n$ algebra we have deduced 
\cite{IsOg3} the Hamiltonians 
\be
\lb{ham5}
{\cal H} = \sum_{m=1}^{n-1} \left( \sigma_m + \frac{(q-q^{-1}) \nu}{\nu +a} \kappa_m \right) + 
\frac{(q-q^{-1}) \xi}{y_1 - \xi}  \; ,
\ee
where $\xi^2 = -a \, c/\nu$ and the parameter $a$ can take one of two values $a=\pm q^{\pm 1}$. Now 
different local representations $\rho$ of the algebra $\alpha BMW_n$ give different integrable spin 
chain models with Hamiltonians $\rho({\cal H})$ which in particular possess $U_q(osp(N|K))$ 
symmetries for some $N$ and $K$. So, representations 
$\rho$ of the algebra $\alpha BMW_n$ are related to the spin chain models of $osp$ - type with $n$ sites 
and nontrivial boundary conditions. BMW chains (chains based on the BMW algebras) describe in a 
unified way spin chains with $U_q(osp(N|K))$ symmetries.

The Hamiltonians for Hecke chain models are obtained from Hamiltonians for BMW chain models by 
taking the quotient $\kappa_j =0$. These models were considered in \cite{IsOg2}, \cite{Isa2}. The 
Hecke chains (chain models based on the Hecke algebras) describe in a unified way spin chains 
with 
$U_q(sl(N|K))$ symmetries. In \cite{IsOg2}, \cite{Isa2} we investigated the integrable open 
chain models formulated in terms of generators of the Hecke algebra (nonaffine case, $y_1=1$). 
  For the open Hecke chains of finite size,
the spectrum of the Hamiltonians with free boundary conditions is determined \cite{IsOg2} for 
special (corner type) irreducible representations of the Hecke algebra. In 
\cite{Isa2} we investigated the functional equations for the transfer matrix type elements of the Hecke 
algebra appeared in the theory of Hecke chains.

We postpone to future publications a construction of the algebra which extends the BMW algebra by
the free algebra with generators labeled by the oscillating Young tableaux (as it is done for the
Hecke algebras in \cite{OgPya}).

\vskip .5cm
\noindent {\bf Acknowledgements.} The work of A. P. Isaev was partially supported by the
grants RFBR 08-0100392-a, RFBR-CNRS 07-02-92166-a and RF President Grant N.Sh. 195.2008.2.


\begin{thebibliography}{99}

\bibitem{Brau}  R. Brauer,
{\it On algebras which are connected with the semisimple continuous groups}, Ann. Math. {\bf 38}
(1937) 854-872.

\bibitem{M1} J. Murakami, {\em The Kauffman polynomial of links and representation theory},
Osaka J. Math. {\bf 24} (1987) 745--758.

\bibitem{BW} J. S. Birman and H. Wenzl, {\em Braids, link polynomials and a new algebra}, Trans.
Amer. Math. Soc. {\bf 313} no. 1 (1989) 249--273.

\bibitem{Ogi} O. Ogievetsky, {\em Uses of Quantum Spaces}, in: Quantum symmetries in theoretical
physics and mathematics, Contemp. Math., {\bf 294} (2002) 161-231; Amer. Math. Soc., Providence,
RI.

\bibitem{FRT}  L. D. Faddeev,  N. Yu. Reshetikhin and L. A. Takhtajan,
{\em Quantization of Lie groups and  Lie algebras}, (Russian) Algebra i Analiz
 {\bf 1}  no.1 (1989) 178--206. English translation in: Leningrad Math. J. 
 {\bf 1}  no.1 (1990) 193--225.

\bibitem{Isa} A. P. Isaev, {\em Quantum groups and Yang-Baxter equations},
preprint MPIM (Bonn), MPI 2004-132,
(http://www.mpim-bonn.mpg.de/html/preprints/preprints.html);
Sov.\ J.\ Part.\ Nucl.\  {\bf 26} (1995) 501
(Fiz. Elem. Chastits i At. Yadra {\bf 26} (1995) 1204).

\bibitem{IOP} A. P. Isaev, O. V. Ogievetsky and P. N. Pyatov, {\em On $R$-matrix representations
of Birman-Murakami-Wenzl algebras}, Trudy Mat. Inst. Steklova {\bf 246} (2004) 147-153. English 
translation in: Proceedings of the Steklov Institute of Mathematics {\bf 246} (2004) 134–141; 
ArXiv: math.QA/0509251.

\bibitem{W2} H. Wenzl, {\em Quantum Groups and Subfactors of Type B, C and D},
Comm. Math. Phys. {\bf 133} (1990) 383-432.

\bibitem{BB} A. Beliakova and Ch. Blanchet, {\em  Skein
construction of idempotents in Birman-Murakami-Wenzl algebras}, Mathematische Annalen
 {\bf 321} no.2 (2001) 347-373; ArXiv: math.QA/0006143.

\bibitem{WT} I. Tuba and H. Wenzl, {\em On braided tensor categories of type BCD},
J. Reine Angew. Math. {\bf 581} (2005) 31-69; ArXiv: math.QA/0301142.

\bibitem{OWen} R. Orellana and H. Wenzl, {\em q-Centralizer Algebras for Spin Groups}.
J. Algebra {\bf 253} no.2 (2002) 237-275.

\bibitem{OV} A. Okounkov and A. Vershik, {\em A new approach to representation theory
of symmetric groups}, Selecta Mathematica, New Ser. {\bf 2} no.4 (1996) 581-605.

\bibitem{IsOg} A. P. Isaev and O. V. Ogievetsky, {\em On representations of Hecke algebras},
Czech. J. of Physics, {\bf 55} no.11 (2005) 1433-1441. \\
 A. P. Isaev and O. V. Ogievetsky, {\em Representations of A-type Hecke algebras},
 Proceedings of International Workshop "Supersymmetries and Quantum Symmetries", Dubna, 2006;
 arXiv: 0912.3701[math.QA].

\bibitem{IsOg5} A. P. Isaev and O. V. Ogievetsky, {\em Representations of Birman-Murakami-Wenzl
algebras}, in preparation.

\bibitem{Nazar} M. Nazarov, {\em Young's orthogonal form for Brauer's Centralizer
Algebra}, Journ. of Algebra {\bf 182} (1996) 664--693.

\bibitem{IsOg3} A. P. Isaev and O. V. Ogievetsky,
{\em On Baxterized Solutions of Reflection Equation and Integrable Chain Models}, Nucl. Phys.
 {\bf B 760} [PM] (2007) 167-183; arXiv: math-ph/0510078.

\bibitem{IsOg2} A. P. Isaev, O. V. Ogievetsky and A. F. Os'kin,
{\em Chain models on Hecke algebra for corner type representations}, Reports on Math. Phys. 
 {\bf 61} no.2 (2008) 309 - 315; arXiv: 0710.0261[math.QA].

\bibitem{Isa2} A. P. Isaev, {\em Functional equations for transfer-matrix operators in open
Hecke chain models}, Theor. Math. Phys.  {\bf 150} no.2 (2007) 187-202.

\bibitem{OgPya} O. Ogievetsky and P. Pyatov, {\em Lecture on Hecke algebras},
in Proc. of the Int. School. "Symmetries and Integrable Systems", Dubna (1999).

\end{thebibliography}
\end{document}